\documentclass[11pt]{article}
\usepackage{graphicx} 
\usepackage{ragged2e}
\usepackage{setspace}
\usepackage{amsmath}
\usepackage{pxfonts}
\usepackage{geometry}
\usepackage{setspace}

\title{Mosaic number and Tile number of Corner Connection Tiles}
\author{Vincent Lin}
\date{September 2023}

\begin{document}

\maketitle

\clearpage

\begin{center}\begin{Large}
    \textbf{Abstract}
\end{Large}\end{center}

\begin{singlespace}    
Lomonaco and Kauffman introduced knot mosaics in 2008 to model physical quantum states. These mosaics use a set of tiles to represent knots on \textit{n} x \textit{n} grids. In 2023 Heap introduced a new set of tiles that can represent knots on a smaller board for small knots. Completing an exhaustive search of all knots or links, \textit{K}, on different board sizes and types is the most common way to determine invariants for knots, such as the smallest board size needed to represent a knot, $m(K)$, and the least number of tiles needed to represent a knot, $t(K)$. In this paper, we propose a solution to an open question by providing a proof that all knots or links can be represented on corner connection mosaics using fewer tiles than traditional mosaics $t_c(K) < t(K)$, where $t_c(K)$ is the smallest number of corner connection tiles needed to represent knot \textit{K}. We also define bounds for corner connection mosaic size, $m_c(K)$, in terms of crossing number, $c(K)$, and simultaneously create a tool called the \textit{Corner Mosaic Complement} that we use to discover a relationship between traditional tiles and corner connection tiles. Finally, we construct an infinite family of links $L_n$ where the corner connection mosaic number $m_c(K)$ is known and provide a tool to analyze the efficiency of corner connection mosaic tiles.
\end{singlespace} 

\clearpage

\section{Introduction}
\begin{figure}[h!]
    \centering
    \includegraphics[width=1\linewidth]{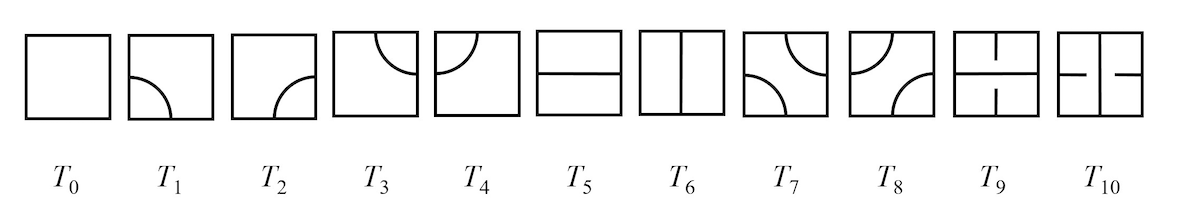}
    \caption{The set of tiles used to construct traditional knot mosaics.}
    \label{fig:1}
\end{figure}

In 2008, Lomonaco and Kauffman \cite{lomonaco2008quantum} introduced the concept of knot mosaics to model physical quantum states in their paper \textit{Quantum knots and mosaics}. Mosaic knot theory uses a set of 11 tiles shown in Figure \ref{fig:1} to create a projection of a knot or link.

\begin{figure}[h!]
    \centering
    \includegraphics[width=1\linewidth]{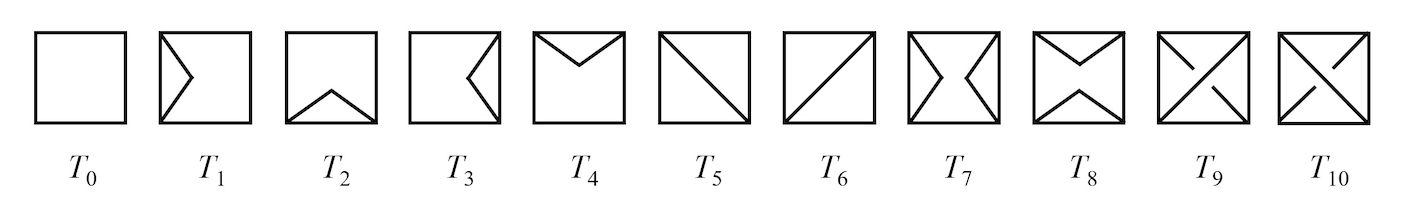}
    \caption{The set of tiles used to construct corner connection knot mosaics.}
    \label{fig:2}
\end{figure}

Lomonaco and Kauffman defined Reidenmeister-like moves (a set of three moves that manipulate a knot without changing the knot type) and conjectured that tame knot theory is equivalent to knot mosaic theory. In other words, two knots are of the same type if and only if there exists a series of Reidenmeister moves relating their mosaic projections. This was later proven by Kuriya and Shehab \cite{kuriya2014lomonaco}. While much work has been done in traditional knot theory, it is reasonable to ask if there are different sets of tiles that could better model knots and provide more powerful invariants; there has been some exploration in hexagonal tiles by Bush \cite{bush2020hexagonal} and Howard \cite{howards2019infinite}. In this paper, we explore a set of tiles introduced by Heap \cite{heap2023knot} shown in Figure \ref{fig:2}.

\begin{figure}[h!]
    \centering
    \includegraphics[width=1\linewidth]{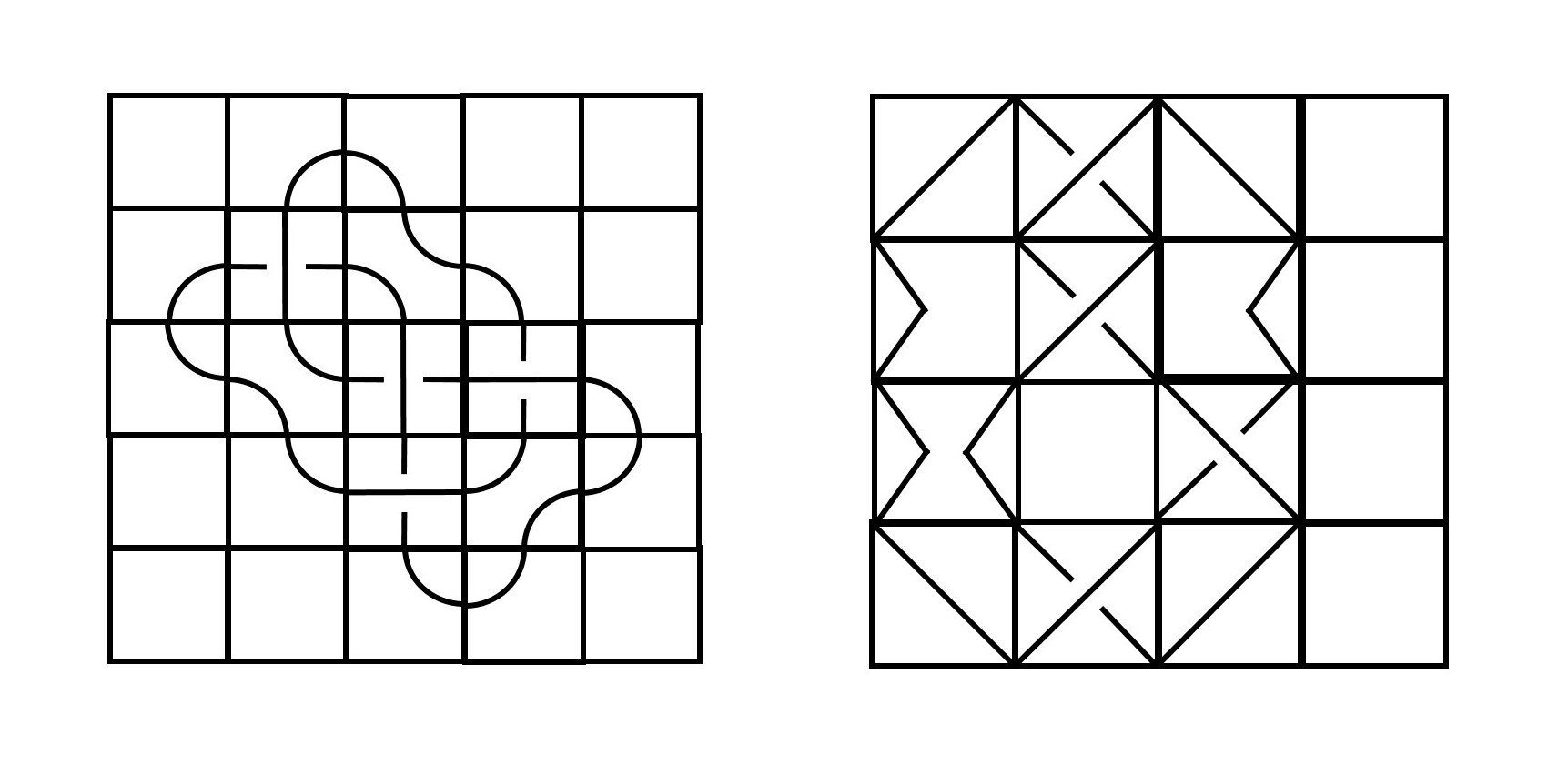}
    \caption{Projection of the $4_1$ knot on traditional and corner connection mosaic, respectively.}
    \label{fig:3}
\end{figure}

The $4_1$ knot can be projected on traditional tiles on a 5 by 5 board using 17 tiles, whereas it can be projected on corner connection mosaics on a 4 by 4 board using only 11 tiles as shown in Figure \ref{fig:3}. In fact, all knots with 8 crossings or less have been tabulated through an exhaustive search on corner connection mosaics, with the result being that all knots 8 crossings or less can be represented on corner connection tiles more efficiently in terms of using less blank tiles \cite{heap2023knot}. This result naturally prompts the question of whether all knots can be represented on corner connection mosaics using fewer tiles.

We claim that for any knot \textit{K}, the smallest number of corner connection tiles needed to represent \textit{K}, $t_c(K)$, is always less than the smallest number of traditional tiles needed to represent \textit{K}, t(K); in other words: 

$$
t_c(K) < t(K).
$$

In section 2, we introduce some preliminary knot theory terminology that we will use in our discussion that follows. In section 3, we will create a tool called \textit{Corner Mosaic Complement} that we will then use to answer an open question from Heap \cite{heap2023knot} in section 4. In section 5, we create bounds for crossing numbers in terms of mosaic number. Finally, in section 6, we introduce a family of links where the corner connection mosaic number is always known.

\section{Preliminaries}
We introduce knot theory terminology, as given by Adams \cite{adams1994knot}:

\noindent \textbf{Definition.} (Knot) A \textit{knot} denoted \textit{K} is a closed curve in 3-space that does not intersect itself anywhere. We do not distinguish between the original closed knotted curve and the deformations of that curve through space that do not allow the curve to pass through itself. The different pictures of the knot that result from these deformations are called \textit{projections} of the knot. 

Invariants are tools used to classify knots. One of the most common ways is the \textit{crossing number}:

\noindent \textbf{Definition.} (Crossing number) The \textit{crossing number} of a knot \textit{K} is the minimal number of crossings in any projections of \textit{K}, denoted \textit{c(K)}.

We can also observe a collection of knots, known as \textit{links}:

\noindent \textbf{Definition} (Links) A \textit{link} is a set of knots in which the knots do not intersect each other but can be tangled together. Each knot that makes up a link is called a \textit{component}.

\noindent \textbf{Definition.} (Split Links) A \textit{split link} is a link whose components can be deformed so that they lie on different sides of a plane in 3-space.

Now, we introduce terminology specific to knot mosaics and corner connection knot mosaics as given by Heap \cite{heap2020space} \cite{heap2023knot}.

\noindent \textbf{Definition.} (Connection Point) We call the midpoint of the edges of traditional tiles, or the corners of corner-connection tiles, a \textit{connection point} if it is also the endpoint of a curve drawn on that tile.

\noindent \textbf{Definition.} (Suitably Connected) A tile in a mosaic is said to be \textit{suitably connected} if all of its connection points touch a connection point on another tile.

\noindent \textbf{Definition.} (\textit{n-mosaic}) An \textit{n} x \textit{n} array of tiles is an \textit{n} x \textit{n knot mosaic}, or \textit{n-mosaic}, if each of its tiles is suitably connected.

\noindent \textbf{Definition.} (Mosaic Number) We define \textit{mosaic number} as the smallest integer \textit{n} such that \textit{K} can fit on a \textit{n}-mosaic using traditional tiles, denoted \textit{m(K)}, or corner connection tiles, denoted \textit{$m_c$(K)}.

\noindent \textbf{Definition.} (Tile Number) We define \textit{tile number} as the smallest number of tiles needed to construct \textit{K} on any size mosaic using traditional tiles, denoted \textit{t(K)}, or corner connection tiles, denoted \textit{$t_c$(K)}.

\noindent \textbf{Definition.} (\textit{k}-Submosaic) We define a sub-mosaic as a \textit{k}-submosaic if it is a submatrix for a \textit{n}-mosaic, where $n \geq k$.

\begin{figure}[h!]
    \centering
    \includegraphics[width=.5\linewidth]{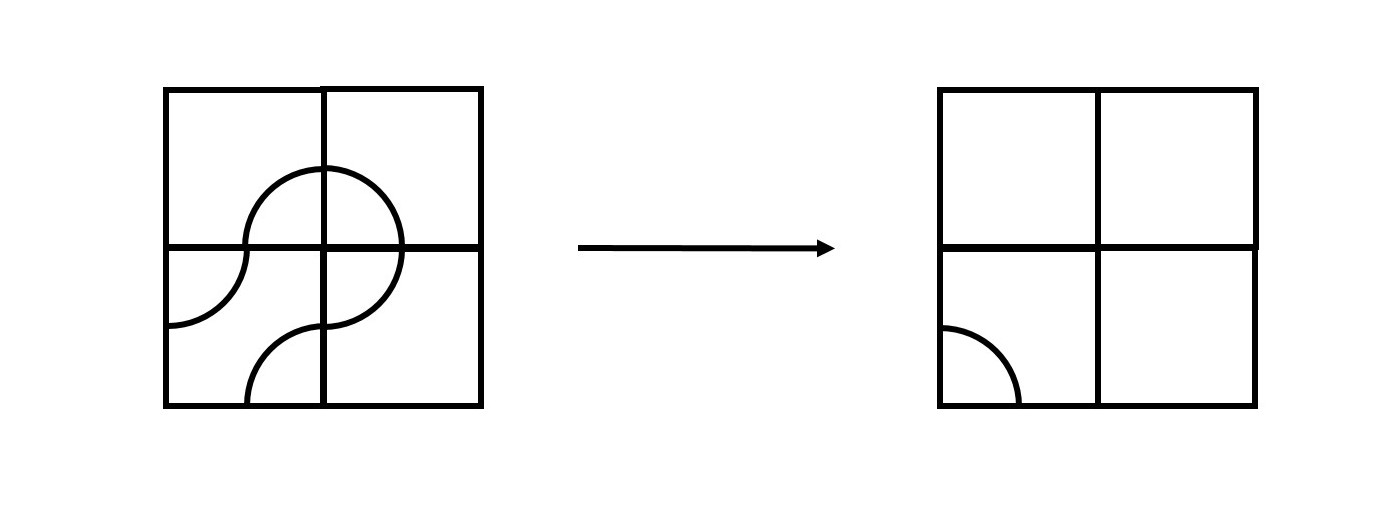}
    \caption{Example of a planar isotopy move.}
    \label{fig:4}
\end{figure}

While working with knot mosaics, we can move knots around via \textit{mosaic planar isotopy moves}. An example of a mosaic planer isotopy move is given in Figure \ref{fig:4}. We can replace any 2 x 2 submosaic with any of the two submosaics without changing the knot types. Throughout this paper, we will be using mosaic planar isotopy moves for traditional tiles and corner connection tiles to use fewer non-blank tiles.

\clearpage

\begin{figure}[h!]
    \centering
    \includegraphics[width=0.75\linewidth]{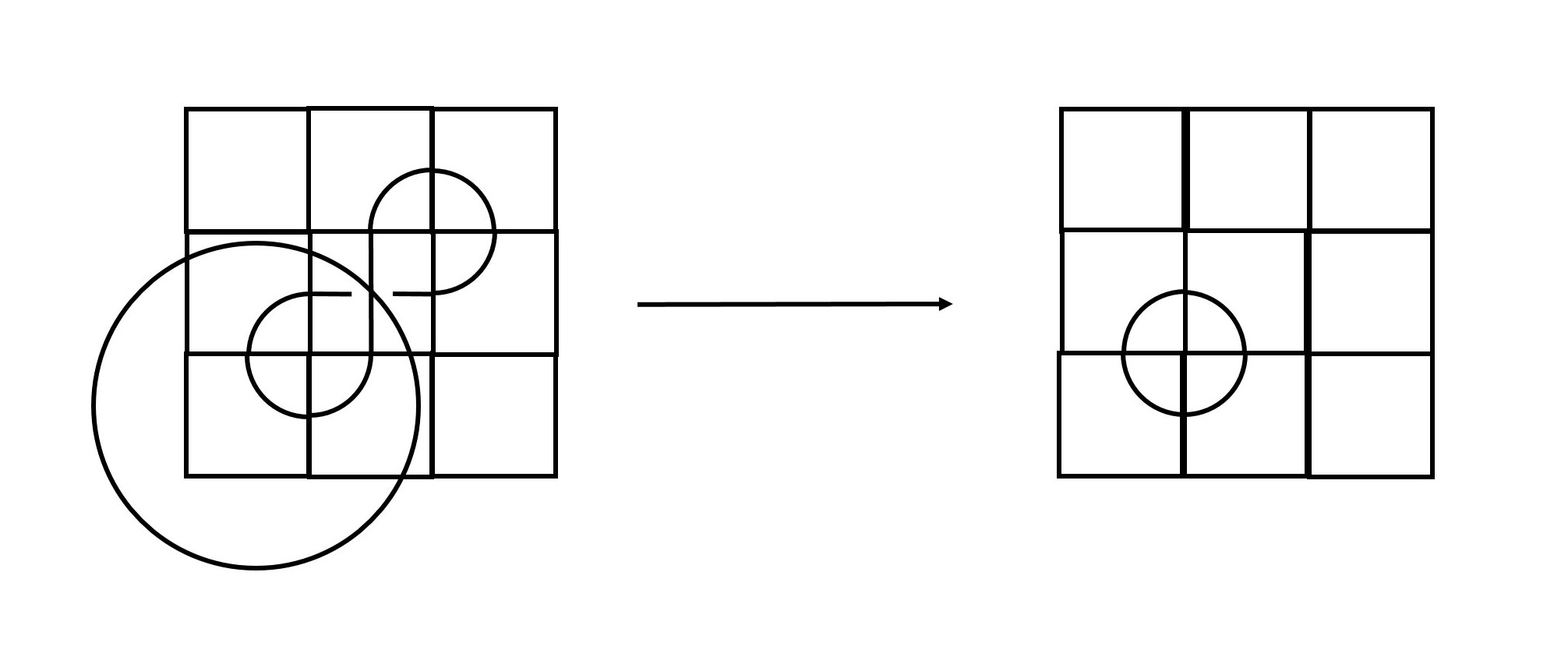}
    \caption{Depiction of a reducible unknot with a circle that meets the diagram transversely at the crossing.}
    \label{fig:5}
\end{figure}

\noindent \textbf{Definition.} (Reducible) A crossing in a knot diagram is \textit{reducible} if there is a circle in the projection plane that meets the diagram transversely at the crossing but does not meet the diagram at any other point as shown in Figure \ref{fig:5}.

\noindent \textbf{Definition.} (Reduced) A knot mosaic is considered \textit{reduced} if there are no reducible crossings on a knot diagram.

\noindent \textbf{Definition.} (Space-efficient) A knot \textit{n}-mosaic is \textit{space-efficient} if it is reduced and the number of non-blank tiles is as small as possible without changing the knot type of the depicted knot.

\noindent \textbf{Notation:} A tile on a mosaic can be denoted \textit{$A_{i,j}$} where \textit{i} is the row of the tile and \textit{j} is the column of the tile.

\section{Construction of a Corner Mosaic Complement} 

The most common way to determine tile number and mosaic number for traditional and corner connection mosaics is to complete an exhaustive search of all possible knots and combinations on a certain \textit{n}-mosaic. We offer a new tool to analyze tile numbers of space-efficient nontrivial knots and non-split links more efficiently. By creating a projection of \textit{K} on corner connection tiles while being equal to the knot type from traditional tiles, we can better analyze the tile number and mosaic number.

To begin our construction, we begin with two results from Heap \cite{heap2018tile} that will assist us with creating the Corner Mosaic Complement. 

\noindent \textbf{Lemma 3.1}. (Heap \cite{heap2018tile}) Suppose we have a space-efficient \textit{n}-mosaic with $n \geq$ 4 and no unknotted, unlinked link components. Then the four corner tiles are blank $T_0$ tiles (or can be made blank via a planar isotopy move that does not change the tile number). The same result holds for the first and last tile location of the first and last occupied row and column.

\begin{figure}[h!]
    \centering
    \includegraphics[width=0.75\linewidth]{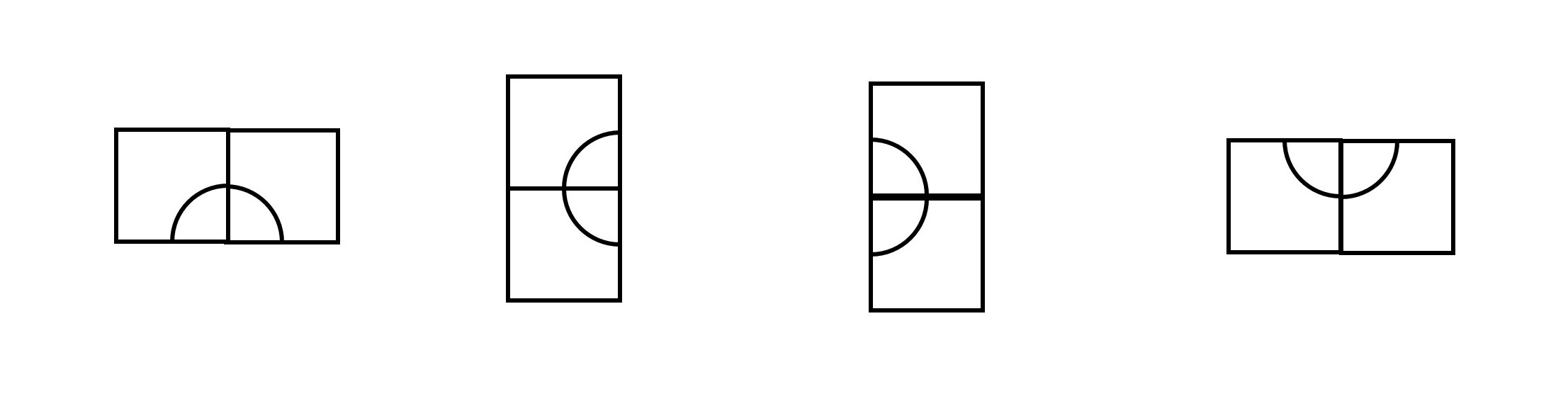}
    \caption{Depiction of tops caps, left caps, right caps, and bottom caps respectively.}
    \label{fig:6}
\end{figure}

\noindent \textbf{Lemma 3.2}. (Heap \cite{heap2018tile}) Suppose we have a space-efficient \textit{n}-mosaic of a knot or link. Then the first occupied row of the mosaic can be simplified so that the non-blank tiles form only top caps. In fact, there will be \textit{k} top caps for some \textit{k} such that 1 $\leq k \leq (n-2)/2$. Similarly, the last occupied row is made up of bottom caps, and the first and last occupied columns are made up of left caps and right caps, respectively. (See figure \ref{fig:6} for caps)

\clearpage

\subsection{Construction of a Corner Mosaic Complement for $\textit{n}\geq5$}

\begin{figure}[h!]
    \centering
    \includegraphics[width=.5\linewidth]{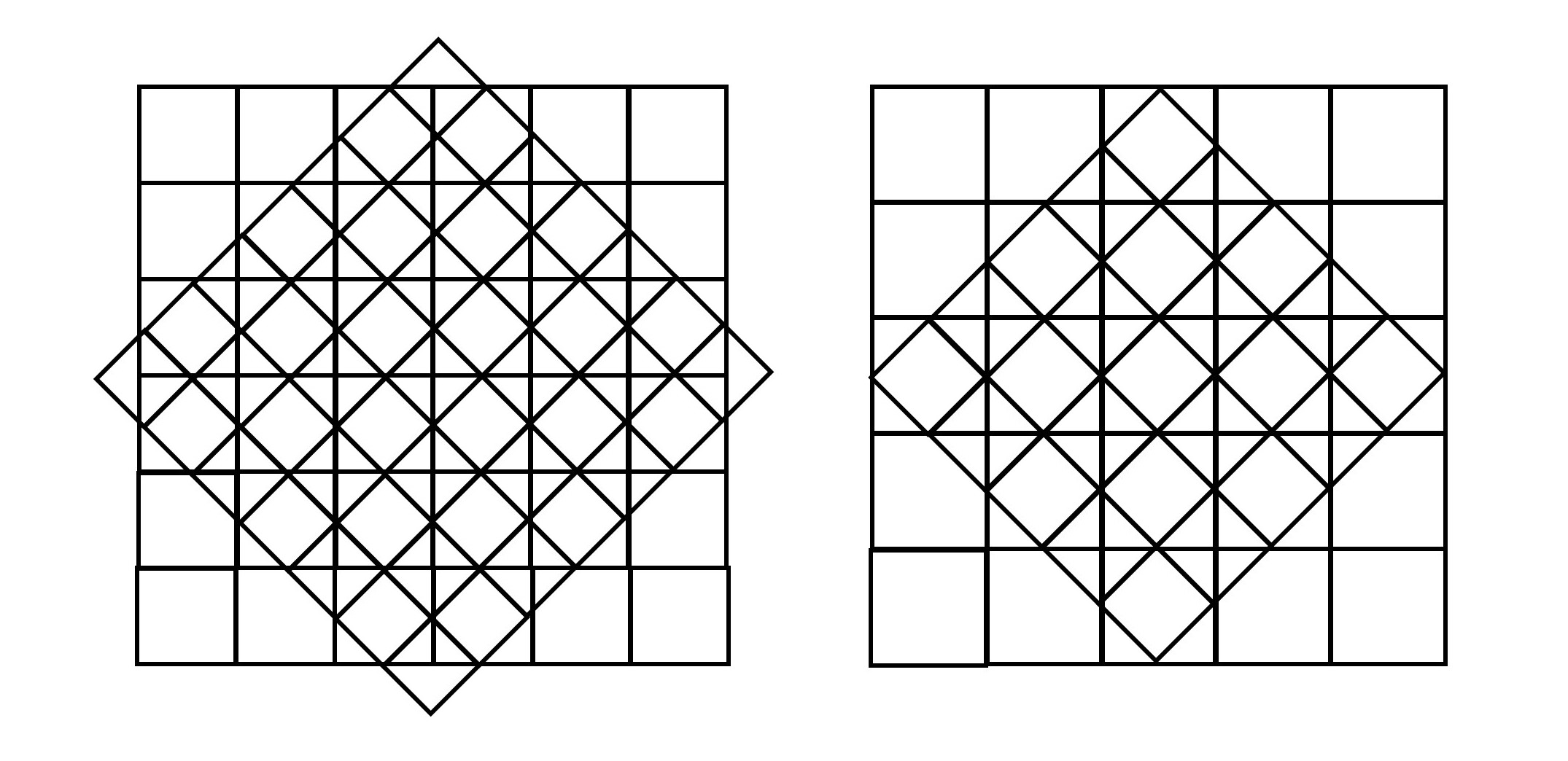}
    \caption{Example of a Corner Mosaic Complement created from a traditional 6-mosaic (left) and a traditional 5-mosaic (right). The size of the Corner Mosaic Complement are 7x7 and 5x5 respectively.}
    \label{fig:7}
\end{figure}

The goal of the Corner Mosaic Complement is to create a Corner Connection Mosaic from traditional mosaics. First, we start with a \textit{n}-mosaic for $\textit{n}\geq5$, place a point at the midpoint of the top edge at $A_{1,3}$ and $A_{1,n-2}$. Similarly, we place points at rotations. Place a point at the midpoint of the right edge for $A_{3,n}$ and $A_{n-2,n}$; a point at the midpoint of the bottom edge for $A_{n,3}$ and $A_{n,n-2}$; and a point at the midpoint of the left edge for $A_{3,1}$ and $A_{n-2,1}$. Four lines are drawn to connect the points to form a square in a 45 degree angle while placing points at intersecting lines. Finally, draw lines through the points on the tilted square to create an array, the \textit{n}-mosaic. Finally, we assume that the tiles for the traditional mosaic are suitably connected, there are no trivial knots or split links, and the knot or link depicted is a projection of the knot that is, based on Lemma 3.2, space-efficient with only top, left, right, and bottom caps. Figure \ref{fig:7} illustrates one example.

\noindent \textbf{Lemma 3.3}. All tiles from the set of traditional tiles have a corresponding tile from the set of corner connection tiles.

\noindent \textbf{Proof.} Take a traditional tile, place a point at the midpoint of each of its edges, and connect them to form an inscribed square. Each traditional tile that isn't $T_0$ has connection points on the midpoints of the edges. If the inscribed square is a representation of a corner connection tile, its connection points would be from the corners of the inscribed square, which are also the midpoints of the edges of the circumscribed traditional tile. As shown in Figure \ref{fig:8}, all tiles from the set of traditional tiles match with the set of tiles from connection corner tiles. 

\begin{figure}[h!]
    \centering
    \includegraphics[width=1\linewidth]{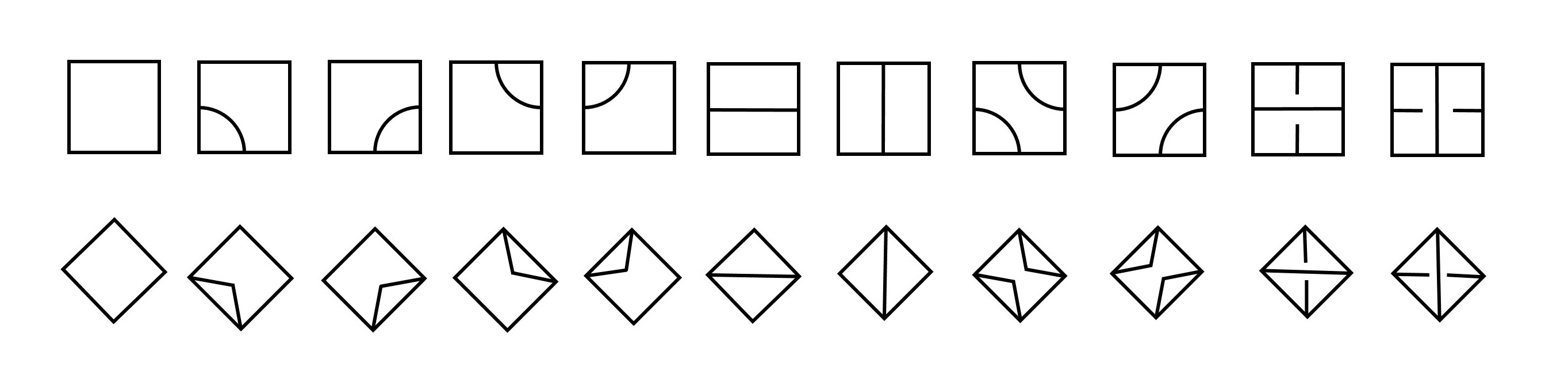}
    \caption{All tiles from the set of traditional tiles and their corresponding tile from the set of corner connection tiles.}
    \label{fig:8.}
\end{figure}

From the construction of the inscribed mosaic as shown in Figure \ref{fig:7}, we can see that most of the tiles, except for the the corner tiles and their adjacent tiles, have an inscribed tile. From Lemma 3.3, we know that each of the tiles with an inscribed tile has a corresponding tile from the set of corner connection tiles. From Lemma 3.1, we can leave the corner tiles without inscribed square because they will be blank $T_0$ tiles. 

\noindent \textbf{Lemma 3.4}. Tiles adjacent to corner tiles do not need inscribed squares. 

\begin{figure}[h!]
    \centering
    \includegraphics[width=1\linewidth]{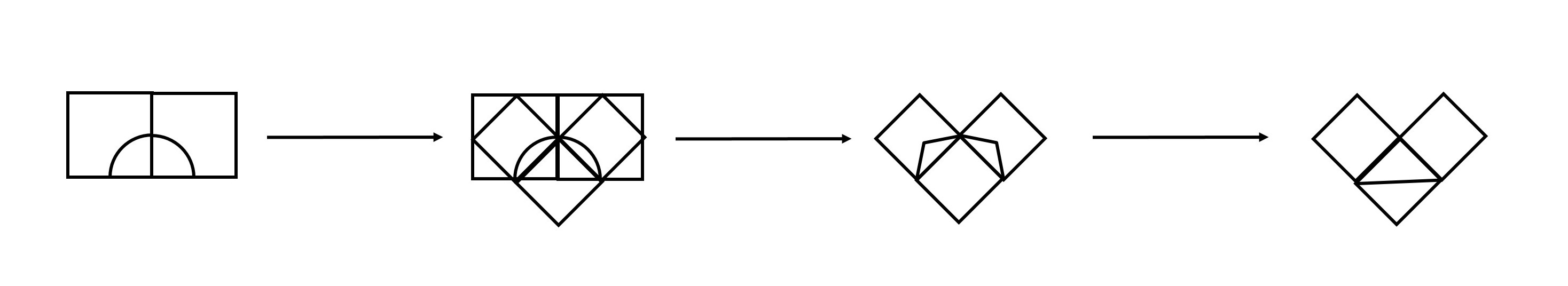}
    \caption{Depiction of caps, their corresponding corner connection tiles, and a planar isotopy move.}
    \label{fig:9}
\end{figure}

\noindent \textbf{Proof.} We know from Lemma 3.2 that the top row can only form top caps. This means that the only tiles possible are $T_0$, $T_1$, and $T_2$ tiles. By Lemma 3.1, the corner tiles must be blank $T_0$ tiles. This leaves only two cases where the tiles adjacent to the corner tiles will have non-blank tiles. As shown in Figure \ref{fig:9}, we can take a cap from traditional tiles and manipulate it into a single non-blank tile in corner connection mosaic. For the first case, there is a top cap to the right of a corner tile in positions $A_{1,2}$ and $A_{1,3}$. From Figure \ref{fig:9}, we can observe that for all caps, the Corner Mosaic Complement could be manipulated through planar isotopy moves to make the inscribed tile at $A_{1,2}$ a blank tile. When creating the Corner Mosaic Complement, we can exclude this tile to form a smaller mosaic as shown in Figure \ref{fig:7}, where the inscribed tiles in tile $A_{1,2}$ is a blank $T_0$ tile. For the second case, there is a top cap to the left of a corner tile in positions $A_{1,n-2}$ and $A_{1,n-1}$. We can apply the same logic for case one and exclude the inscribed tile in $A_{1,n-1}$. In fact, we can apply this logic to all caps through rotation; by rotating the mosaic 90 degrees clockwise, we can make the right caps become top caps and apply the same logic to exclude the inscribed tiles and repeat for bottom and left caps. Finally, if the tiles adjacent to the corner tiles are blank tiles, then we can exclude their corresponding corner connection tiles, since they will be blank $T_0$ corner connection tiles. 

Recall in Section 3.1, we discussed the construction of a Corner Mosaic Complement for $n\geq5$; we now consider the cases where $n\leq4$.

\noindent \textbf{Lemma 3.5}. The Corner Mosaic Complement for a traditional \textit{n}-mosaic where $n\leq3$ does not exist, and when \textit{n} = 4, we have a 3-mosaic. 

\begin{figure} [h]
    \centering
    \includegraphics[width=0.5\linewidth]{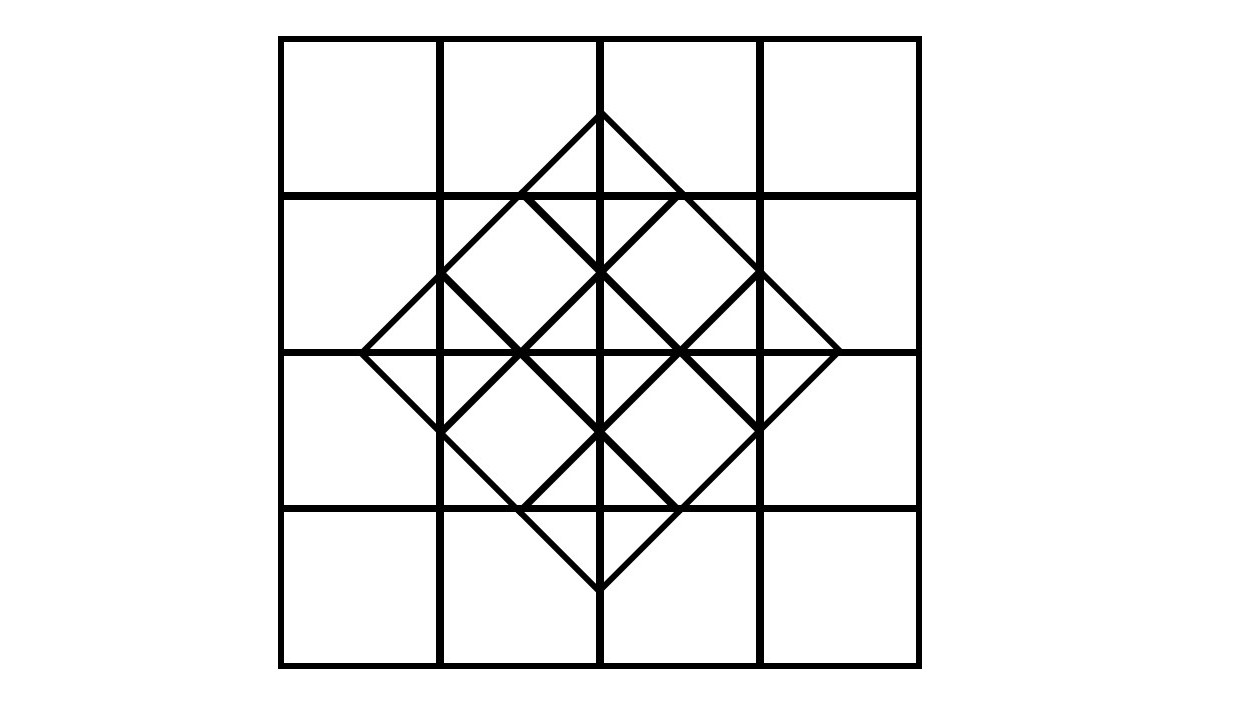}
    \caption{Depiction of the Corner Mosaic Complement for a 4-mosaic.}
    \label{fig:10}
\end{figure}

\noindent \textbf{Proof.} We know that there does not exist a projection of a non-trivial knot or non-split link that can be depicted on traditional mosaics for $n\leq3$ \cite{lomonaco2008quantum}. As we are assuming no trivial knots and no non-split links, we therefore do not need to construct a Corner Mosaic Complement for traditional \textit{n}-mosaics where $n\leq3$. Now consider a 4-mosaic. We can create inscribed squares on the inner four tiles and then create a final mosaic resulting in a  3 x 3 square by abiding by the rules outlined in Lemmas 3.1, 3.3, and 3.4, as shown in Figure \ref{fig:10}.

\noindent \textbf{Lemma 3.6}. Knots or links depicted by the Corner Mosaic Complement are of the same knot or link type as the original knot projected on a traditional mosaic.

\noindent \textbf{Proof.} Consider again Figure \ref{fig:7}; we can observe that connection points for each curve in each tile are at the same spots and are suitably connected to other tiles in the same way. For caps, the resulting curve after placing it on the blank corner connection has the same connection points as connection points of the caps. Since the curve doesn't pass through itself, we know that the caps are also suitably connected in the same way. We can generalize this for any \textit{n}-mosaic, as all non-trivial non-split link knots can be projected on a traditional mosaic with caps and tiles not on the perimeter of the mosaic. Therefore this allows all tiles not on the perimeter to have inscribed tiles suitably connected to other inscribed tiles in the same way and with the caps being placed on blank corner connection tiles.

This completes the construction for Corner Mosaic Complement.

\section{Corner Connection Tile Number}

In this section, we propose a proof for the open question from Heap\cite{heap2023knot}:

\noindent \textbf{Question}. Is it always true that $t_c(K) \leq t(K)$?

\noindent \textbf{Theorem 4.1}. For all knots and non-split links, the corner connection tile number is less than the traditional mosaic tile number, $t_c(K) < t(K)$.

\begin{figure} [h]
    \centering
    \includegraphics[width=0.5\linewidth]{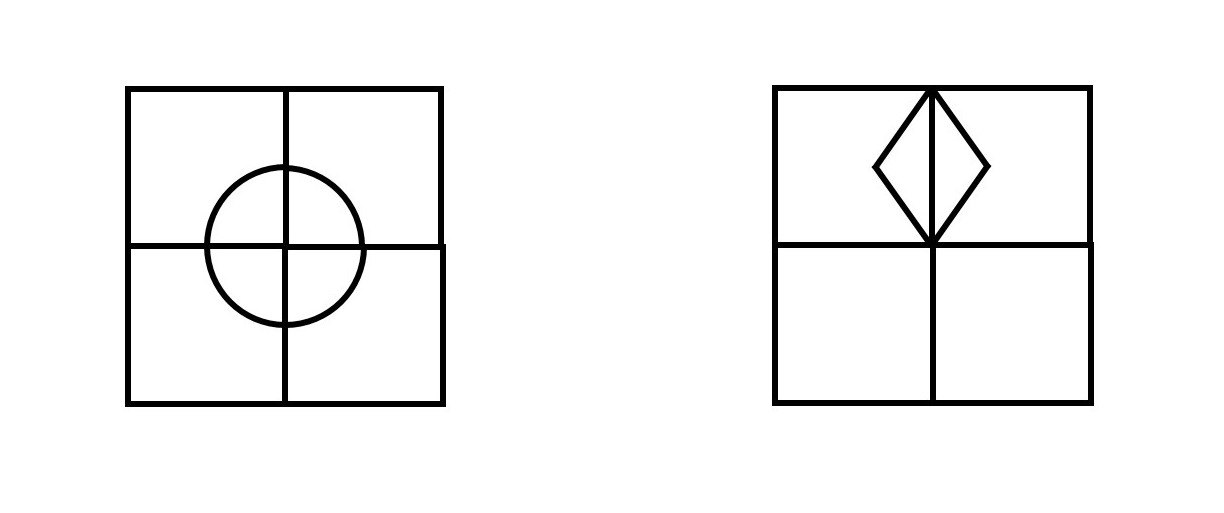}
    \caption{Projection of the unknot on a traditional mosaic with its mosaic number and tile number realized on the left. Projection of the unknot on a Corner Connection Mosaic with its mosaic number and tile number realized on the right.}
    \label{fig:11}
\end{figure}

\noindent \textbf{Proof}. From Lemma 3.6, we know that a space-efficient knot depicted on a traditional mosaic is equivalent to its Corner Mosaic Complement. By Lemma 3.2, we know that there exists a projection of a knot or non-split link on a traditional mosaic where there are only caps on the first and last rows and columns. By Lemma 3.3 other tiles of the knot can be represented by a corner connection tile. By Lemma 3.4 we can place caps that use 2 traditional tiles on one tile from the set of corner connection tiles. Since there exists a space-efficient tile for every knot or non-split link with caps on the first row on traditional mosaics, the Corner Mosaic Complement can always be created with fewer non-blank tiles. (For knots on mosaic sizes $n\leq3$, the only knot that has a projection is the unknot on 2-mosaic and 3-mosaic. However its tile number is 4, and it can be represented on corner-connection tiles with just two tiles, as shown in Figure \ref{fig:11}.)

\clearpage

\noindent \textbf{Remark}. Although there exists a projection of any non-trivial knots or non-split links on Corner Connection Mosaic that uses fewer tiles than a space-efficient tile, the Corner Connection Mosaic itself isn't always space-efficient or shown on the smallest mosaic size. For example, the $4_1$ knot can only be projected on a traditional 5-mosaic, therefore its Corner Mosaic Complement would be a 5-mosaic. However, as shown in Figure \ref{fig:3}, the mosaic number of the $4_1$ knot, $m_c(4_1)$, is 4.

Finally, we introduce a new tool to prove that split-links can also be projected on a corner connection mosaic using fewer tiles.

\subsection{Construction of an Inefficient Corner Mosaic Complement}

We can create an inefficient corner mosaic complement by placing inscribed squares inside every tile within a traditional mosaic. We then create a square formed in a 45 degree angle that includes the inscribed tiles. In other words, the inefficient Corner Mosaic Complement can be made using the process of creating a Corner Mosaic Complement instead of making the corner connection mosaic smaller using Lemma 3.1, Lemma 3.3, and Lemma 3.4. 

\noindent \textbf{Theorem 4.2} For all split-links, the corner connection tile number is less than the traditional mosaic tile number, or $t_c(K) < t(K)$.

\noindent \textbf{Proof.} We note that all tiles from the traditional mosaic will have a corresponding mosaic tile complement, thus all projections of split links on a traditional mosaic can be projected onto the Inefficient Corner Mosaic Complement. The Inefficient Corner Mosaic Complement of the split-link will be suitably connected with the same logic as the proof for Lemma 3.6 and be of the same knot type for each of its link components. We can conclude that is $t_c(K) \leq t(K)$. We can sharpen this relationship by using Lemma 3.2 to realize that we can reduce the link component's caps to use one fewer tile using the planar isotopy move described in Figure \ref{fig:9}. By Lemma 3.2, there must exists caps that can be reduced using planar isotopy moves described in Figure \ref{fig:9}, thus always resulting in fewer tiles than its projection on traditional tiles. Thus, $t_c(K)<t(K)$.

\section{Bounds for Corner Connection Tiles}

It is challenging to identify bounds for the possible crossing number of knots that can be projected on any given \textit{n}-mosaic. We must prove that the knots with crossing number lower than the lower bound must have a projection on a smaller mosaic, and knots with crossing number higher than the upper bound cannot be projected on the mosaic. Previous work on creating a lower bound on the crossing number in terms of mosaic number utilized a system called the \textit{grid diagram}\cite{lee2014mosaic}. This paper will utilize this bound as it is crucial to the construction of a bound for corner connection tiles. We first state Theorem 5.1, proven by Lee et. al. \cite{lee2014mosaic}.

\noindent \textbf{Theorem 5.1} \cite{lee2014mosaic}. Let \textit{K} be a nontrivial knot or a non-split link except the Hopf link, then $m(K) \leq c(K) + 1$. 

In Theorem 5.2 and 5.3, we introduce a new naming convention, where \textit{n} refers to the \textit{n}-mosaic created from traditional tiles, and $n_c$ refers to the \textit{n}-mosaic created from corner connection tiles. We also introduce a definition \textit{inner tiles} that is used in the proof of Theorem 5.2.

\noindent \textbf{Definition} (Inner tiles) The \textit{inner tiles} are defined as tiles that are not in the perimeter of the mosaic. More formally, they are tiles that have 8 tiles total directly adjacent and diagonally from their positions in the mosaic.

To begin creating bounds for corner connection tiles, we first observe that the invariant crossing number of a knot, \textit{c(K)}, remains constant for traditional tiles and corner connection tiles. We begin with a theorem that defines a relationship between mosaic number $m(K)$ and corner connection mosaic number $m_c(K)$.

\noindent \textbf{Theorem 5.2}. For all space-efficient projections of knots and non-split links \textit{K} on a \textit{n}-mosaic where $n\geq4$, there exists a projection of \textit{K} on a $n_c$-mosaic where:

$$
n_c \leq 2n-5.
$$

\noindent \textbf{Proof}. We note that in creating a Corner Mosaic Complement, the inner tiles will always have an inscribed tile. Since there is not an inscribed tile on the corners of the \textit{n}-mosaic, the corners of the inscribed squares from the set of inner tiles are on the perimeter of the Corner Mosaic Complement. The size of the corner mosaic complement is 1 less than the sum of
the number of rows and number of columns. Finally, we know that the corner connection mosaic number may not be realized on the Corner Mosaic Complement, so the corner connection mosaic number may be less than the size of the Corner Mosaic Complement.

\indent Thus we have
\begin{equation*}
    (n-2)+(n-2)-1=2n-5.
\end{equation*}

\noindent \textbf{Theorem 5.3}. The upper bound of $m_c(K)$ can be bounded by the crossing number by:

$$
m_c(K)\leq2c(K)-3.
$$

\noindent \textbf{Proof}. By Theorem 5.1, we know that for traditional tiles, $m(K)$ is bounded above by $c(K)$. We know that from Theorem 5.2, the upper bound of the $n_c$-mosaic needed to project a knot \textit{K} on corner connection \textit{n}-mosaic is 2\textit{n}-5. We can use the upper bound of $m(K)$ with respect to the crossing number and the upper bound needed to represent knot K on a traditional mosaic on a traditional board to find the the upper bound of $m_c(K)$ in terms of the crossing number.

We then have:
\begin{align*}
m(K) & \leq c(K)+1, \\
m_c(K) & \leq 2(c(K)+1)-5, \\
m_c(K) & \leq 2c(K)-3.
\end{align*}

\section{Infinite Family of Links where the Mosaic Number is Known}

It is known that the upper bound for the crossing number in terms of mosaic number grows faster than the upper bound for crossing number in terms of corner mosaic number \cite{heap2023knot}. In other words, for sufficiently large knots, there exists a knot or link that has a mosaic number less than the corner mosaic number. Since the rate at which the bounds for traditional tiles compared to corner connection mosaics grows much faster, it may seem intuitive that as the crossing number of knot or link \textit{K} grows large, it can be projected on a smaller traditional mosaic simply because traditional mosaics can contain more crossing tiles than corner connection mosaics of equal size when n is large. In fact, Heap \cite{heap2023knot} proved that some large knots have mosaic number less than corner mosaic number $m(K) \leq m_c(K)$. In this section of this paper, we will construct a family of links and describe its special properties to provide more tools to answer a question we propose.

\noindent \textbf{Question 6.1}. Does there exist an infinite number of knots or links \textit{K} where the $m_c(K) \leq m(K)$?

We begin by creating the link defined below.

\noindent \textbf{Definition 6.2}. We define $L_n$ as an alternating link on a $n_c$ - mosaic where $n_c = n$  and $n_c = 2k+1,$ for $k \in \mathbb{Z}$ with crossing tiles in positions

\begin{align*}
\{(i,j): 1 \leq i \leq n\} & \backslash \{(2i,2j): i \leq \lfloor n/2 \rfloor, j \leq \lfloor n/2 \rfloor \} \\
& \backslash \{(2i-1,2j-1): i \leq \lceil n/2 \rceil, j \leq \lceil n/2 \rceil \}.
\end{align*}

\begin{figure}[h]
    \centering
    \includegraphics[width=0.40\linewidth]{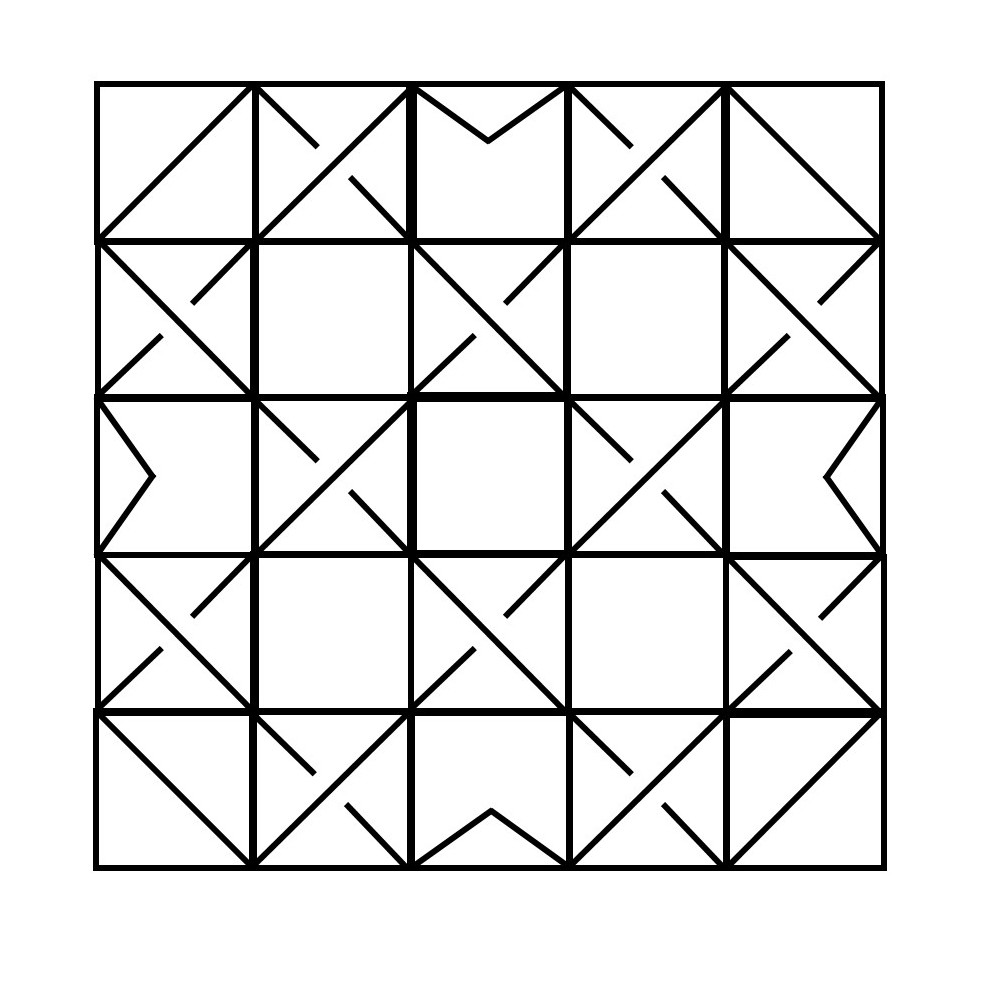}
    \caption{Example of a $L_5$ link.}
    \label{fig:12}
\end{figure}

In other words, $L_n$ is projected in a chain-like pattern on a corner connection as shown in Figure \ref{fig:12}.

We now establish properties of this infinite family of links by first introducing the famous Thistlethwaite theorem and a lemma about Corner Connection tiles.

\noindent \textbf{Theorem 6.3}. (Kauffman\cite{kauffman1988new}, Thistlethwaite\cite{thistlethwaite1987spanning}, Murasugi\cite{murasugi1987jones}) Any reduced diagram of a link has minimal crossings.

\noindent \textbf{Lemma 6.4}. (Heap \cite{heap2023knot}) For any $n \geq 3$, the upper bound for the number of crossing tiles used in an $n_c$-mosaic created from corner connection tiles is $n^2/2$ if \textit{n} is even and $(n^2+n-4)/2$ if \textit{n} is odd.

\noindent \textbf{Theorem 6.5}. The crossing number for $L_n$ is the number of $T_9$ and $T_{10}$ tiles, or simply, Link $L_n$ is reduced, with crossing number $c(L_n) = \lfloor n^2/2 \rfloor$

\noindent \textbf{Proof}. We can create the projection of $L_n$ with alternate crossings by placing $T_9$ tiles every odd row and $T_{10}$ tiles every even row as the the crossing tiles. By Thistlethwaite's Theorem \cite{kauffman1988new}, this link would be a reduced projection. Therefore the number of crossing tiles is equivalent to $L_n$'s crossing number $c(K)$ with the construction of $L_n$ as defined in definition 6.2 always containing $\lfloor n^2/2 \rfloor$ crossing tiles.

\noindent \textbf{Theorem 6.6}. Link $L_n$ has $\lceil n/2 \rceil$ link components.

\noindent \textbf{Proof}. Each link component in a corner connection mosaic representation of $L_n$ will have 2 possible projections as shown in Figure \ref{fig:13}:

\begin{figure} [h]
    \centering
    \includegraphics[width=0.75\linewidth]{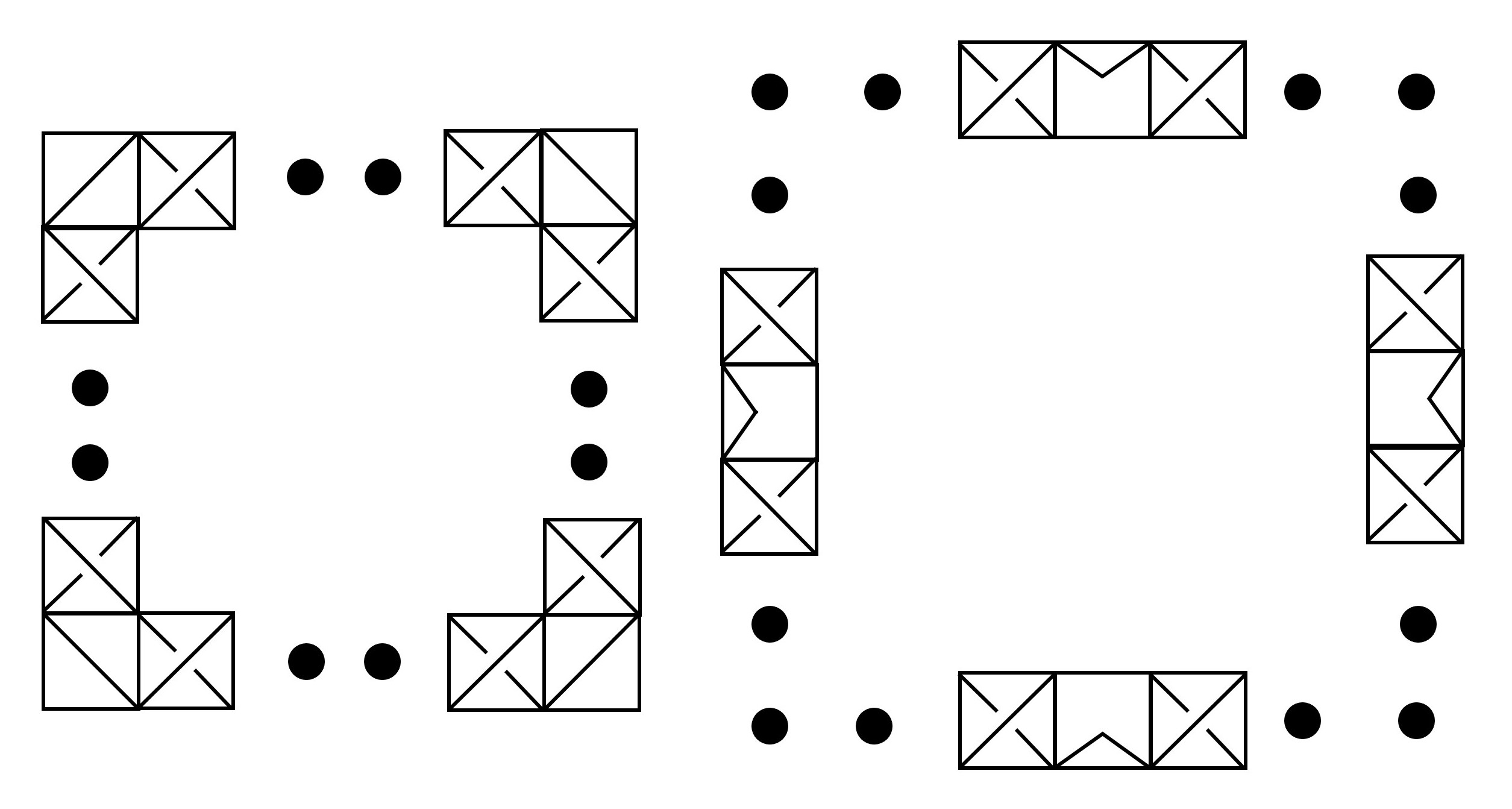}
    \caption{The two projections for each link component of $L_n$.}
    \label{fig:13}
\end{figure}

If we observe only the first column of the left mosaic in figure 13, all link components will have a $T_1$, $T_5$, or $T_6$ tile. In the right mosaic of figure 13, every even tile in the first column is either $T_1$, $T_5$, $T_6$ or a crossing tile. We know from the construction of $L_n$ that every other tile in the first column is a crossing tile, where the first tile is $T_6$. We can therefore count the number of other tiles in the first column to count the number of link components, which is $\lceil n/2 \rceil$.

\noindent \textbf{Theorem 6.7}. Link $L_n$ has mosaic number \textit{n}, or $m(L_n)=n$.

\noindent \textbf{Proof}. To begin our proof by contradiction, suppose $L_n$ could be projected on a $(n-1)_c$-mosaic. The $(n-1)_c$-mosaic would be even since $L_n$ is originally projected on a odd \textit{n} by \textit{n} mosaic. The maximum number of crossing tiles $(n-1)_c$-mosaic can fit can be found using Theorem 6.4

\[
(n-1)^2/2.
\]

Link $L_n$ will always have more crossing tiles than the maximum number of tiles possible one mosaic smaller than the it, $(n-1)_c$-mosaic when comparing the crossing number from Theorem 6.5 to the maximum number of crossings tiles possible. Because $L_n$ can be projected on an $n_c$-mosaic as shown when constructing $L_n$. We have reached a contradiction.

$L_n$ can be always be projected as a reduced knot on traditional mosaics as shown in Figure 14. Simultaneously it can also be projected on corner connection mosaics as a smaller mosaic than the one we present in Figure 14. Because there are no obvious space-efficient planar isotopy moves, we propose the following conjecture:

\begin{figure} [h]
    \centering
    \includegraphics[width=0.5\linewidth]{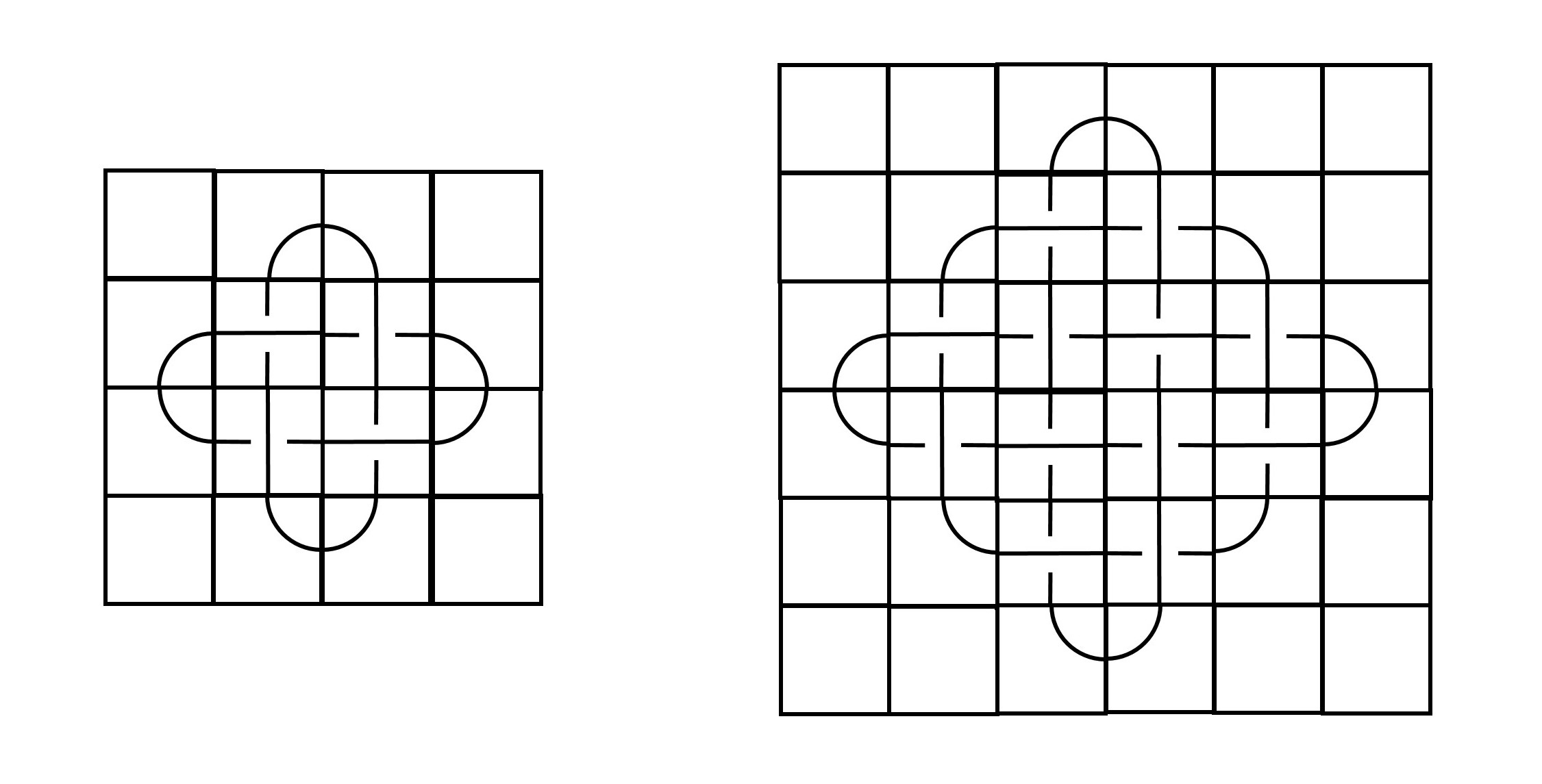}
    \caption{Example of $L_3$ (left) and $L_5$ (right).}
    \label{fig:14}
\end{figure}

\noindent \textbf{Conjecture 6.8}. The corner connection connection mosaic of $L_n$ is less than the mosaic number of $L_n$. $m_c(L_n) \leq m(L_n)$. 

\section{Future Work}

We develop a tool to analyze corner connection mosaic efficiency and generalized the corner connection tile number of knots in comparison to traditional tiles. Future work can try to improve on the bounds and theorems proposed in this paper. 

In particular, space-efficient knots on traditional mosaics have multiple caps, and each cap can be manipulated to reduce the tile by one. Can we create a better relationship between the corner connection tile number and traditional tile number using this idea?

Further work can also be done to improve on the bounds proposed by this paper. None of the knots tabulated currently have corner connection mosaic number at its upper bound. We can investigate ways to create a stricter upper bound or find a relationship between corner connection mosaic number and crossing number with respect to mosaic number.

Finally, we can continue to improve our understanding of the invariants corner connection mosaics produce, tabulation of knots on corner connection mosaics allows us to know the mosaic number and tile number of knots, as well as interesting properties the mosaic projections that has these invariants realized may show.

\clearpage

\section{Summary}
Knot theory has applications in many different fields. For example, we can use knots to understand the behavior of knotted DNA and its relation to topoisomerase, an enzyme with crucial roles in DNA replication and transcription, to create chemo drugs to combat cancer such as doxorubicin. Knots can also be used to study chirality and isomers within molecular structures, as well as to control the stability of a molecule\cite{kruve2019ion}. For example, there has been previous research in creating knotted molecules to possess certain properties. As shown by Lomonaco and Kauffman, we can also use knots to model physical quantum states \cite{lomonaco2008quantum}. Invariants are important in knot theory to distinguish between multiple knots and their properties. Knot mosaics are useful as they introduce a new set of invariants such as tile number and mosaic number, as well as the ability to study knots by representing them on mosaics and  matrices. Finding elementary proofs and creating tools such as Corner Mosaic Complement to analyze knot mosaics and their properties without using an exhaustive search is therefore useful to compute invariants of knots in more general cases.

\clearpage

\section{Acknowledgements}
I am deeply grateful to Dr. Avineri and Dr. Boltz for their proofreading and invaluable suggestions to this paper. Their attention to detail and dedication to enhancing the clarity and precision of this work have been instrumental in its overall quality.

I am truly grateful to Dr. Bullard for supporting my curiosity in my first computer science project. The experience from that project gave me the research skills that made this project possible.

\bibliographystyle{plain}

\bibliography{mybibfile}
All images in this paper were created by the author of this paper.

\end{document}